# Nonhomogeneous Subordinators and Lévy Copulas for Reliability Analysis of Time-Dependent Degradation Processes

Yu Shi[1, *], Qianmei Feng[2]

[1]*Department of Mathematics, Tarleton State University, Stephenville, TX 76401*

[2]*Department of Industrial and Systems Engineering, University of Houston, Houston, TX 77204*

## ABSTRACT

Many systems in engineering applications may experience complex multi-dimensional degradation processes involving random jumps. Analyzing the reliability of such degradation processes presents significant challenges. Previously developed models have utilized multi-dimensional Lévy subordinators and Lévy copulas to describe these processes, with reliability functions expressed through characteristic functions and Laplace transforms (Shu *et al*., 2015, 2016, Shi *et al*., 2020, Shi and Feng, 2024). However, these models are restricted to stationary linear paths and homogeneous Lévy subordinators. In real-world scenarios, increments in a degradation process are often time-dependent and nonhomogeneous. To address this limitation of previous studies, this work extends existing degradation models by introducing nonhomogeneous Lévy subordinators, which are designed to capture time-varying complexities inherent in degradation processes. Unlike conventional degradation models that often assume a constant rate of jumps, nonhomogeneous Lévy subordinators aim to provide a more flexible representation of degradation behaviors that evolve under varying external conditions and internal system dynamics. Additionally, the joint behaviors of multiple degradation processes are characterized using Lévy copulas, enabling a more comprehensive analysis of system reliability. Numerical examples are provided to illustrate the applicability of these models in lifetime analysis.

**Corresponding author:** Yu Shi, Department of Mathematics, Tarleton State University, Stephenville, TX 76402, USA.

**Email:** YSHI@tarleton.edu

**Phone:** 254-968-1766

## 1. Introduction

Degradation processes in various engineering systems, such as railways, aircraft, oil and gas equipment, and civil infrastructure, evolve randomly due to a combination of internal properties and varying external conditions (MacLean *et al*., 2018, Shi *et al.*, 2020). These random events may result in failures, causing significant production and safety issues as well as financial losses for organizations. In more complicated cases, degradation processes often exhibit time-dependent random jumps, resulting in non-stationary increments and non-linear behaviors. Such behaviors are frequently observed in complex systems where the degradation rate changes over time due to internal and external factors. For instance, the growth rate of internal cracks in civil infrastructure can be influenced by time-dependent factors (Boyajian, 2002). To address these complexities, it is necessary to study non-decreasing degradation processes with time-dependent increments by extending the time-homogeneous models to nonhomogeneous ones.

Stochastic processes offer a powerful mathematical framework for capturing inherent uncertainties and variabilities associated with the degradation of engineering systems (Singpurwalla, 1995, Liu *et al*., 2020). Since the 1970s, stochastic process modeling has been extensively studied for degradation and reliability analysis. The Wiener process is one of the first stochastic models applied to degradation modeling. Initially used to describe stock price volatility in the financial sector, the Wiener process also forms the basis of the Black-Scholes pricing formula (Black and Scholes, 1973). In engineering contexts, the Wiener process is an effective tool for modeling non-monotonic changes in degradation processes. To address its limitations in representing cumulative degradation with monotonically non-decreasing paths, stochastic processes such as Gamma and inverse Gaussian processes have been introduced (Liu *et al*., 2020). Gamma processes, for example, are characterized by independent increments that follow a Gamma distribution.

While these stochastic processes effectively model degradation with independent increments, many complex random jump mechanisms in degradation processes remain inadequately described by such approaches. In response to this limitation, Lévy processes have been introduced to model degradation processes involving random jumps. Pruitt (1981) first demonstrated that Lévy processes provide key insights into the asymptotic behavior of random walks. Lévy processes are a type of stochastic processes characterized by their independent and stationary increments, which enable them to play a fundamental role in modeling various phenomena that involve sudden jumps (Mozumder, 2011). Shu *et al*. (2015) explored the reliability function and lifetime characteristics of non-decreasing degradation processes with random jumps, utilizing Lévy measures to capture various jump mechanisms in these processes. Recent advancements have extended these foundational ideas. Shi *et al.* (2020) constructed general multi-dimensional Lévy processes governed by Lévy copulas, effectively handling uncertainties and dependencies among degradation processes. Their work provides robust reliability predictions for multi-dimensional degradation systems. However, the degradation processes discussed thus far are limited to time-homogeneous scenarios.

The degradation processes discussed thus far are limited to time-homogeneous scenarios, which do not account for time-dependent random jumps. In real-world applications, however, it is common to encounter systems exhibiting multiple non-stationary degradation processes. For instance, stock and asset prices are often time-dependent and uncertain, making them suitable for modeling with nonhomogeneous subordinators. While such processes are not strictly Lévy subordinators due to the lack of stationary increments, they can be regarded as a generalization of subordinators (Cont and Tankov, 2004, Orsinger *et al.,* 2016). For example, nonhomogeneous Poisson processes have been utilized to analyze dependent degradation processes, providing analytical expressions for system reliability and proposing dynamic inspection policies based on system health information (Castro and Landesa, 2024). Nonhomogeneous subordinators extend the concept of homogeneous Lévy subordinators by relaxing the assumption of stationary increments previously studied (Abdel-Hameed, 2014). Unlike homogeneous Lévy subordinators, nonhomogeneous subordinators are defined by a set of time-dependent Lévy measures $v_t$ (Beghin and Ricciuti, 2016). This new Lévy measure is a function of time, dynamically capturing jump behaviors at each time point. The primary advantage of nonhomogeneous subordinators is their ability to extend stationary increments to non-stationary conditions, enabling the modeling of more complex and general degradation processes. This approach offers greater flexibility compared to homogeneous Lévy subordinators.

Building upon existing stochastic modeling approaches using Lévy processes, this paper extends current methodologies by introducing a novel framework to capture multi-dimensional, time-dependent degradation processes. This is achieved by utilizing multi-dimensional Lévy copulas and nonhomogeneous subordinators. The proposed approach relaxes the limitations of traditional models, which often assume independence or homogeneity of degradation increments. The new framework provides a more flexible and realistic approach that can capture complex dependencies and time-varying behaviors. Furthermore, this study offers a significant advancement in reliability analysis by deriving closed-form expressions for the reliability function, providing practical tools for assessing the performance and predicting the lifetime of complex systems.

This paper is organized as follows. Section 2 introduces key concepts and preliminaries related to additive processes and nonhomogeneous subordinators. In Section 3, we present multi-dimensional Lévy copulas and nonhomogeneous subordinators to characterize the joint behaviors of multi-dimensional, time-dependent degradation processes. Section 4 derives the reliability function. Numerical examples illustrating the proposed models are provided in Section 5, and conclusions are presented in Section 6.

The notation used in this paper is listed in Table 1.

Table 1: List of Notation

| | |
|---|---|
| $X_t$ | An additive process |
| $X_{St}$ | A nonhomogeneous subordinator |
| $m$ | Dimension of an additive process / nonhomogeneous subordinator |
| $X_{St}^i$ | The $i^{\text{th}}$ dimension of a nonhomogeneous subordinator |
| $v_s(dy)$ | Dynamic Lévy measure at time $s$ |
| $\psi(t,\lambda)$ | Laplace exponent of a nonhomogeneous subordinator |
| $\varphi_{X_{St}}(u)$ | Characteristic function of an $m$-dimensional nonhomogeneous subordinator |
| $\eta_s(u)$ | Lévy-Khintchine representation for $\forall s \in [0, t]$ |
| $F_{X_{St}}(x)$ | Cumulative distribution function (CDF) |
| $f_{X_{St}}(x)$ | Probability density function (PDF) |
| $T_x$ | First passage time |
| $R(t)$ | Reliability function |
| $C_t$ | Copula function at time step $t$ |
| $U(x)$ | Tail integral on Lévy measure |
| $\alpha$-$MS(t)$ | $\alpha$-multi-stable process |
| $\theta$ | The dependence parameter in Clayton copula |

## 2. Preliminaries of Nonhomogeneous Subordinators

Additive processes are generalizations of Lévy processes, achieved by relaxing the stationary condition (Cont and Tankov, 2004, Abdel-Hameed, 2014). As a subclass of additive processes, nonhomogeneous subordinators are characterized by non-decreasing sample paths, which have gained attention in recent years (Orsinger *et al.,* 2016, Beghin and Ricciuti, 2016). In this section, we introduce key preliminaries of additive processes and nonhomogeneous subordinators, focusing on their formulation in $R^m$ space.

**Definition 1** *A stochastic process* $X_t$ *in* $R^m$ *space is called an additive process if it is cadlag* $X_0 = 0$, *and satisfies the following properties (Cont and Tankov, 2004):*

1. *Independent increments: for every increasing sequence of times* $t_0, \ldots, t_n$, *the random variables* $X_{t1} - X_{t0}, \ldots, X_{tn} - X_{t,n-1}$ *are independent.*

2. *Stochastic continuity:* $\forall \varepsilon > 0, \forall\, t > 0, \lim_{h \to 0} P(|X_{t+h} - X_t| \geq \varepsilon) = 0$.

With the relaxation of time stationarity, the measure of jumps must be redefined to accommodate time dependence. This nonhomogeneous measure can also be converted into a homogeneous one by applying a continuous deterministic time transformation. Let $Y(s)$ be a Lévy process in $R^m$ space such that $s$ is a time variable $s \in [0, t]$, and $v_t$ be a continuous function such that $v_t(0) = 0$ and $\int_0^\infty (y \wedge 1)\, v_t(dy) < \infty$. Then the additive process $X_t$ can be rewritten as

$$X_t(s) = Y_{v_t}(s). \tag{1}$$

Then, the characteristic function of $X_t$ at a fixed time $t$ is (Kluge, 2005):

$$E[e^{iuX_t}] = e^{\int_0^t \eta_s(u) ds} \tag{2}$$

where $\eta_s(u) = -\frac{1}{2}u^2 A_s + iub(s) + \int_{R^+} \left[e^{iuy} - 1 - iuy 1_{|y| \leq 1}\right] v_s(dy)$ is the dynamical Lévy-Khintchine representation at time step $s \in [0, t]$. $A_s$ is a positive definite $m \times m$ matrix, $b(s)$ is a constant drift term, and $v_s(dy)$ is the dynamic Lévy measure at time step $s$.

Similar to regular Lévy processes in $R^m$, an additive process also possesses an infinitely divisible distribution at any time $t$, uniquely determined by *time spot characteristics* $\left(A_t, v_t, b(t)\right)_{t\geq 0}$ (Cont and Tankov, 2004). For a family of time-dependent dynamical Lévy measure $v_s$ on time interval [0, $t$], it satisfies $\int_0^t \int (1 \wedge |y|^2)\, v_s(dy)ds < \infty$ (Beghin and Ricciuti, 2018). When $A_t = 0$ and $b(t) \geq 0$ at time $t$, $X_t$ simplifies to a nonhomogeneous subordinator $X_{St}$ with a non-decreasing path, making it well-suited for modeling time-dependent cumulative degradation processes with random jump events.

**Theorem 1 (Cont and Tankov, 2004, Sato, 1999)** *For an additive process $X_t$ in $R^m$ space, $X_t$ has an infinitely divisible distribution at any t. The law of $X_t$ is uniquely determined by its time spot characteristics $\left(A_t, v_t, b(t)\right)_{t\geq 0}$ and expressed by Eq.* (2).

Similar to the Lévy-Itô decomposition for a Lévy process, the local jump characteristic of an additive process $X_t$ is described by a time-dependent Poisson random measure on time interval [0, $t$] (Cont and Tankov, 2004, Sato, 1999). The Lévy-Itô decomposition expresses an additive process as a sum of a continuous part and a compensated sum of jumps (Sato, 1999). The Lévy-Itô decomposition for an additive process $X_t$ is given as Theorem 2.

**Theorem 2 (Cont and Tankov, 2004, Sato, 1999)** *Let $X_s$ be an additive process in $R^m$ satisfying $\int_{|y \leq 1|} |y|\, v_s(dy) < \infty$ for all s >0Then, there exists a probability space $(\Omega, \mathbb{F}, P)$, $\Omega \in \mathrm{F}$ with $P[\Omega] = 1$, such that the positions and sizes of jumps of $X_s$ on* [0, $t$] *are described by a Poisson random measure $J_X$ on [0, t] × $R^m$,*

$$J_X = \sum_{s\in[0,t]} \delta(s, X_s), \tag{3}$$

*Where $\delta(s, X_s)$ denotes the Dirac measure (point mass) at location $(s, y) \in [0],[t] \times R^m$, with (time in-homogeneous) intensity given by $v_s(dx)ds$ such that*

$$E[J_X([t_1, t_2] \times A) = \mu_J([t_1, t_2] \times A) = \int_{t_1}^{t_2} v_s(A)ds. \tag{4}$$

Therefore, the nonhomogeneous subordinator can be defined as follows:

**Definition 2 (nonhomogeneous subordinator)** *A nonhomogeneous subordinator $X_{St}$ is characterized by its Laplace exponent $\psi(t,\lambda)$, such that for any $\lambda \geq 0$:*

$$\mathbb{E}\left[e^{-\lambda St}\right] = e^{\psi(t,\lambda)}.$$

*The Laplace exponent can be expressed in terms of the time-dependent Lévy measure $v_t dx$ such that: $\psi(t,\lambda) = \int_0^t \int_0^\infty \left(1 - e^{-\lambda x}\right) v(s,dx) ds$.*

## 3. Distributional Characterization and Reliability Analysis of Multi-dimensional Nonhomogeneous Subordinators

This section develops the distributional framework for multi-dimensional, time-dependent degradation processes modeled by nonhomogeneous subordinators. Let $X_{St} = \{X_{St}^1, \ldots, X_{St}^m\}$ denote an $m$-dimensional nonhomogeneous subordinator, where each component $X_{St}^i$ represents the cumulative degradation of the $i$th component or degradation indicator. Each marginal process is characterized locally by a time-dependent Lévy measure, while dependence among component-wise jumps is introduced through a Lévy copula.

The analysis proceeds from local jump behavior to finite-time system reliability. Section 3.1 first constructs the joint time-dependent Lévy measure from the marginal Lévy measures using a Lévy copula. Section 3.2 then uses the resulting joint Lévy measure to characterize the finite-time joint and marginal distributions through their Lévy–Khintchine representations and characteristic functions. Finally, Section 3.3 connects these finite-time degradation distributions to first-passage reliability for a series degradation system.

Throughout this section, $s$ denotes the local or integration time used to describe the instantaneous jump characteristics, whereas $t$ denotes the observation horizon over which degradation accumulates. This distinction is important because the local Lévy measure may vary with $s$, while the distribution of $X_{S_t}$ reflects the accumulation of these local characteristics over the entire interval [0, $t$].

### 3.1. Dynamic Lévy Copulas and Joint Dynamic Lévy Measures

The concept of copulas and Sklar's theorem provides a general framework for separating marginal behavior from dependence in multivariate distributions (Aas *et al*., 2009). In an ordinary distributional setting, a copula links marginal cumulative distribution functions to their joint distribution. For multidimensional Lévy processes, dependence among component-wise jumps can instead be characterized through a Lévy copula, which links marginal Lévy tail integrals to a joint Lévy tail integral while preserving the marginal Lévy measures. Our previous work applied Lévy-copula-based dependence modeling to multidimensional degradation and reliability analysis under time-homogeneous settings (Shi *et al*., 2020, 2024). In the present study, this framework is extended to nonhomogeneous degradation processes by allowing the marginal Lévy measures to vary over time and, in the general formulation, by allowing the Lévy-copula dependence structure to vary over time as well.

For the $i$th component, let $v_{i,s}(dy_i)$, $i = 1, \dots, m$, denote the marginal dynamic Lévy measure at time $s$. For the non-decreasing degradation processes considered in this study, whose jumps are nonnegative, define the corresponding marginal Lévy tail integral as:

$$U_{i,s}(x_i) = v_{i,s}([x_i, \infty)),\ x_i > 0. \tag{5}$$

The quantity $U_{i,s}(x_i)$ represents the local intensity of jumps in the $i$th degradation component whose magnitudes are greater than or equal to $x_i$ at time $s$.

Let $v_s(dy) = v_s(dy_1, \dots, dy_m)$ denote the joint dynamic Lévy measure of $X_{S_s}$ at time $s$. Its corresponding joint tail integral is defined by

$$U_s(x_1, \dots, x_m) \quad = \nu_s\big([x_1, \infty) \times \cdots \times [x_m, \infty)\big). \tag{6}$$

Thus, $U_s(x_1, \dots, x_m)$ represents the local intensity of joint jumps whose magnitudes are greater than or equal to the specified threshold $x_1, \dots, x_m$.

To accommodate possible time variation in cross-component jump dependence, we consider a time-indexed family of positive Lévy copulas. At each fixed time $s$, dependence among component-wise jumps is described by a m-dimensional Lévy copula $C_s$. The family may vary with $s$; the time-invariant specification is retained as an important special case.

**Definition 3 (Time-dependent Lévy-copula family)** *For each* $s \in [0, t]$*, let* $C_s$ *be a m-dimensional positive Lévy copula* (Kallsen and Tankov, 2006). Thus, $C_s$ *is grounded, m- increasing, and has identity one-dimensional margins. Specifically,* $C_s(u_1, u_2, \dots, u_m) = 0$ *if any* $u_i = 0$ *and* $C_s(\infty, \dots, u_i, \dots, \infty) = u_i, i = 1,2, \dots, m$. *The collection* $\{C_s : s \in [0, t], C_s : [0, \infty]^m \to [0, \infty]\}$ *is called a time-dependent Lévy-copula family if* $s \to C_s(u_1, u_2, \dots, u_m)$ *is measurable for every* $(u_1, u_2, \dots, u_m) \in [0, \infty)^m$.

At each fixed time $s$, the small jumps in a multi-dimensional nonhomogeneous system are estimated by the time-dependent Lévy measure $\nu_s(dy)$ and its marginal $\nu_{i,s}(dy_i)$ in an observation interval [0, $t$]. Applying Sklar's theorem for Lévy copulas pointwise at each time $s$, the joint local tail integral $U_s$ and the marginal local tail integrals satisfy (Shi *et al.*, 2020, 2024)

$$U_s(x1, x2, \dots, xm) \ = C_s(U_{1,s}(x_1), \dots, U_{m,s}(x_m)). \tag{7}$$

Equation (7) extends the Lévy-copula representation to a time-dependent setting. At each time $s$, $C_s$ preserves the marginal time-dependent Lévy measures and specifies the dependence structure among the jumps of the degradation components.

To obtain an explicit expression for the joint local Lévy density, suppose that each marginal local Lévy measure $\nu_{i,s}(dy_i)$ is absolutely continuous and has a density $f_{i,s}(y_i)$, for $i$ = 1, 2,..., $m$.

Assume further that $C_s$ is $m$-times continuously differentiable. By taking the mixed derivative of the tail-integral relationship in Eq. (7) and applying the chain rule, the joint local Lévy measure $v_s(dy)$ is also continuous. Its joint local Lévy density $f_s(y_1, \ldots, y_m)$ is given by:

$$f_s(y_1, \ldots, y_m) = \left. \frac{\partial^m C_s(U_1, U_2, \ldots, U_m)}{\partial u_1 \partial u_2 \ldots \partial u_m} \right|_{u_i = U_{i,s}(y_i)} \prod_{i=1}^{m} f_{i,s}(y_i). \quad (8)$$

Here, $f_s(y_1, \ldots, y_m)$ denotes the joint local Lévy density at time $s$, and $f_{i,s}(y_i)$ denotes the marginal local Lévy density of the $i$th component. We next consider the special case in which the dependence structure is time-invariant; that is, $C_s = C$ for every $s \geq 0$. Let $c$ denote the Lévy copula density, defined as the $m$th-order mixed partial derivative of $C$. Although the Lévy copula $C$ does not vary with time, the joint local Lévy density generally remains time-dependent because the marginal tail integrals and marginal Lévy densities depend on time $s$. Therefore, the cumulative joint Lévy measure over the observation interval $[0, t]$, denoted by $v^c{}_t$, is obtained by integrating the joint local Lévy density with respect to $s$:

$$v_t^c\,(d\boldsymbol{y}) = \left[\int_0^t c\left(U_{1,s}(y_1), \ldots, U_{m,s}(y_m)\right) \prod_{i=1}^{m} f_{i,s}(y_i) ds\right] d\boldsymbol{y}, \quad (9)$$

where

$$c(u_1, \ldots, u_m) = \frac{\partial^m C(u_1, \ldots, u_m)}{\partial u_1 \ldots \partial u_m}$$

is the Lévy copula density.

Thus, even when the Lévy copula is time-invariant, the cumulative joint Lévy measure remains time-dependent because the marginal tail integrals and marginal local Lévy densities vary over time. Eq. (9) aggregates the local jump behavior from time zero to time t and provides the joint Lévy measure used to derive the finite-time distribution in Section 3.2.

### 3.2. Finite-Time Distribution Induced by Time-Dependent Lévy Characteristics

The local Lévy characteristics developed in Section 3.1 describe the instantaneous jump behavior of the degradation process. Reliability analysis, however, requires the finite-time distribution of the accumulated degradation at an observation horizon $t$. Therefore, this subsection establishes the mapping from time-dependent marginalLévy characteristics to the finite-time joint and marginal distributions of the nonhomogeneous subordinator.

Consider an $m$-dimensional degradation process at observation time $t$, $X_{St} = \{X_{St}^1, \ldots, X_{St}^m\}$. Let $b_t$ denote the cumulative drift vector over [0, $t$], and $\nu_t^C$ denote the cumulative joint Lévy measure obtained by integrating the local joint Lévy characteristics over the same interval. Under the time-invariant Lévy-copula specification considered in Eq. (9), the resulting cumulative joint Lévy measure is denoted by $\nu^c{}_t$. Characteristic-function methods have been widely used to derive probability distributions and reliability measures for stochastic degradation processes (Shu *et al.*, 2015, 2016; Shi *et al.*, 2020, 2024). Accordingly, based on the Lévy–Khintchine representation for additive processes, the joint characteristic function of $X_{St}$ is given by:

$$\varphi_{X_{St}}(u) = E[e^{i\langle u, X_{St}\rangle}] = e^{\int_0^t i\langle b_s, u\rangle + \int_{R^m}\left(e^{i\langle u,y\rangle} - 1 - 1_{|y|<1} i\langle u,y\rangle\right)v_s(d\boldsymbol{y})\, ds} \tag{10}$$

where $\eta_s(u)$ denotes the $m$-dimensional time-dependent Lévy–Khintchine exponent at time $s$, determined by the local joint Lévy measures $v_s(d\boldsymbol{y})$.

Eq. (10) determines the finite-time joint distribution of the multidimensional nonhomogeneous degradation process. To obtain the marginal degradation distributions required for reliability analysis, we next consider the characteristic function of each individual component. For the $k$th marginal nonhomogeneous subordinator, the characteristic function at time $t$ is:

$$\varphi_{X_{St}^k}(u) = E[e^{i\left\langle u, X_{St}^k\right\rangle}] = e^{\int_0^t \eta_S^k(u) ds} \quad, k = 1,2.,,,,m \tag{11}$$

where

$$\eta_S^k(u) = ib_S^k u + \int_0^\infty (e^{iuy} - 1)\, v_{k,s}(d\boldsymbol{y}) \tag{12}$$

is the corresponding one-dimensional time-dependent Lévy–Khintchine exponent.

To recover the marginal probability distribution from the characteristic function, we apply Fourier inversion. Shephard (1991) established the relationship between characteristic functions and probability distributions through Fourier inversion in the one-dimensional setting. Following this approach, but replacing the stationary characteristic function with the finite-time characteristic function of the nonhomogeneous subordinator, the marginal probability density function of the $k$th degradation component can be written as:

$$f_{X_{St}^k}(x) = \frac{1}{2\pi}\int_{-\infty}^{\infty} e^{-iux}\varphi_{X_{St}^k}(u)du = \frac{1}{2\pi}\int_{-\infty}^{\infty} \lim_{h\to 0}\frac{sin(hu)}{hu} e^{-iux}\varphi_{X_{St}^k}(u)du. \tag{13}$$

The corresponding marginal cumulative distribution function is

$$F_{X_{St}^k}(x) = \int_0^{x_k} f_{X_{St}^k}(z,t)\, dz \tag{14}$$

where $\varphi_{X_{St}^k}(u)$ is the characteristic function at time $t$. The marginal CDFs in Eq. (14) provide the finite-time degradation probabilities required for the first-passage reliability formulation developed in Section 3.3.

### 3.3. Reliability of Multi-dimensional Nonhomogeneous Degradation Processes

For the $k$th marginal degradation process $X_{St}^k$, we define the first passage time to the prescribed degradation threshold $x_k$ as:

$$T_k = \inf\{t \geq 0 : X_{St}^k > x_k\}. \tag{15}$$

Since $X_{St}^k$ is non-decreasing, the event $\{T_k > t\}$ is equivalent to the event that the cumulative degradation has not exceeded its threshold $x_k$ by time $t$. Therefore, the marginal reliability function of the $k$th component is

$$R_k(t) = P(T_k > t) = P(X_{St}^k \leq x_k) = F_{X_{St}^k}(x_k) \tag{16}$$

where $F_{X_{St}^k}(x_k)$ is the marginal cumulative distribution function derived in Section 3.2. The monotonicity of the nonhomogeneous subordinator is essential in Eq. (16), because it converts the first-passage event into a finite-time distributional event.

For an $m$-dimensional degradation system, let $x = (x_1, \dots, x_m)$ denote the vector of degradation thresholds. In this study, the degradation system is assumed to operate as a series system. Therefore, system failure occurs when any one of the degradation components exceeds its corresponding threshold. The system lifetime is consequently

$$T = min(T_1, T_2, \dots, T_m). \tag{17}$$

Hence, the system reliability function at the observation time $t$ is

$$R(t) = P(T_1 > t, T_2 > t, \dots, T_m > t) = P(X_{St}^1 \leq x_1, \dots, X_{St}^m \leq x_m) = F_{X_{St}}(x), \tag{18}$$

where $F_{X_{St}}(x)$ denotes the finite-time joint cumulative distribution function of the multi-dimensional nonhomogeneous degradation process.

According to Sklar's theorem, the finite-time joint cumulative distribution function can be expressed in terms of its marginal cumulative distribution functions and a distributional copula $C_F$. Therefore,

$$F_{X_{St}}(x) = C_F(F_{X_{St}^1}(x_1), F_{X_{St}^2}(x_2), \dots, F_{X_{St}^m}(x_m)). \tag{19}$$

Based on Eq. (16), Eq. (19) becomes

$$R(t) = C_F\big(R_1(t), \dots, R_m(t)\big). \tag{20}$$

Eq. (20) shows the multi-dimensional system reliability is determined by the marginal reliability functions together with the dependence structure among the degradation components.

It is important to distinguish the distributional copula $C_F$ in Eq. (19) from the Lévy copula introduced in Section 3.1. The Lévy copula operates at the level of local jump measures and determines the dependence structure of componentwise jumps, whereas $C_F$ operates at the level of finite-time degradation distributions. These two copulas therefore represent dependence at different levels of the model.

Equations (16)-(20) show that the multi-dimensional system reliability is determined by two components: the time-dependent marginal degradation distributions and the dependence structure among the finite-time degradation variables. The former is derived from the nonhomogeneous subordinators in Section 3.2, while the latter is represented through the distributional copula $C_F$.

## 4. Special Case: $\alpha$-Multistable Subordinators with Clayton Dependence

The general framework developed in Section 3 applies to multi-dimensional, nonhomogeneous degradation processes with general marginal time-dependent Lévy measures and dependence structures. In this section, we specialize the framework to a tractable class of models in which the marginal degradation processes are represented by $\alpha$-multistable subordinators, and the dependence structure is characterized by the Clayton family.

This specialization separates the two principal components of the model: the time-varying $\alpha$-multistable Lévy measures characterize the non-stationary marginal degradation behavior, while the Clayton dependence specification governs the cross-component dependence. The resulting model provides an analytically tractable setting for deriving marginal and system reliability functions and for investigating the effects of non-homogeneity and dependence on degradation behavior.

### 4.1. $\alpha$-Multistable Marginal Degradation Model

In this study, each marginal degradation process is modeled by a one-dimensional $\alpha$-multistable subordinator. The $\alpha$-multistable subordinator extends the positive $\alpha$-stable subordinator used in the homogeneous setting by allowing the stability parameter to vary over time (Shi *et al.*, 2020). In contrast to the homogeneous case, the stability parameter is allowed to vary with time, resulting in non-stationary increments and time-dependent jump behavior.

**Definition 4 ($\alpha$-multistable subordinator)** For the $k$th degradation component, $X_{St}^{k}$ is an $\alpha$-multistable subordinator if it is a one-dimensional nonhomogeneous subordinator whose local Lévy measure at time $s$ is given by:

$$v_{k,s}(dy_k)=\frac{\alpha_k(s)}{\Gamma(1-\alpha_k(s)){y_k}^{\alpha_k(s)+1}}dy_k,\ \text{where } 0<\alpha_k(s)<1. \tag{21}$$

Here, $\alpha_k(s)$ denotes the time-dependent stability parameter of the $k$th marginal degradation process. Because $\alpha_k(s)$ varies with time $s$, the corresponding local Lévy measure $v_{k,s}(dy_k)$ is time-dependent, resulting in nonstationary degradation increments. If $\alpha_k(s)$ is constant over time, the model reduces to the homogeneous positive $\alpha$-stable subordinator considered in earlier work.

Once the marginal time-dependent Lévy measure is specified, the finite-time characteristic function and marginal degradation distribution can be obtained using the general formulations developed in Section 3.2.

### 4.2. Clayton Lévy-Copula Specification

Having specified the time-dependent marginal degradation processes in Section 4.1, we next characterize the dependence among their component-wise jumps. Following the Lévy-copula framework developed in Section 3.1, a positive Clayton Lévy copula is adopted to link the marginal Lévy tail integrals, thereby describing the dependence structure among simultaneous degradation jumps.

For an $m$-dimensional system, let $U_k \in [0, \infty], k = 1,2,..,m$. The $m$-dimensional positive Clayton Lévy copula is expressed as:

$$C_L(U_1, \dots., U_m, \theta_L) = \left(\sum_{k=1}^{m} U_k^{-\theta_L}\right)^{\frac{-1}{\theta_L}}, \theta_L > 0, \quad (22)$$

where $\theta_L$ is the dependence parameter. Stronger positive dependence among component-wise degradation jumps implies a higher likelihood of simultaneous large degradation increments across components.

In this study, the dependence structure is assumed to be time-invariant, while the marginal Lévy measures are allowed to vary over time. Therefore, $C_s = C_L = C$ for all $s$. The nonhomogeneity of the degradation process is introduced through the time-varying marginal Lévy measures $v_{k,s}(dy_k)$ specified in Eq. (21), whereas the dependence among component-wise jumps is governed by the constant parameter $\theta_L$.

The resulting joint local Lévy measure is derived in Section 4.3 by combining this Clayton Lévy-copula specification with the $\alpha$-multistable marginal Lévy measures defined in Section 4.1.

### 4.3. Joint Local and Cumulative Lévy Measures under Clayton Dependence

As established in Section 3.1, the joint local Lévy measure is obtained by applying the Lévy copula to the marginal Lévy tail integrals. For notational convenience, define the mixed derivative of the Clayton Lévy copula as

$$c_L(u_1, u_2, \dots, u_m; \theta_L) = \frac{\partial^m C_L(u_1, u_2, \dots, u_m, \theta_L)}{\partial u_1, \dots, \partial u_m}. \quad (23)$$

Under the Clayton Lévy-copula specification in Eq. (22), the general expression in Eq. (8) becomes:

$$v_{C,s}(dy) = c_L(U_{1,s}(y_1), \dots, U_{m,s}(y_m); \theta_L) \prod_{k=1}^{m} v_{k,s}(dy_k), \quad (24)$$

where $U_{k,s}(y_k) = v_{k,s}([y_k, \infty))$ denotes the marginal Lévy tail integral associated with the $k$th marginal local Lévy measure $v_{k,s}(dy_k)$ at time $s$.

Eq. (24) separates two distinct features of the degradation model. The time-varying $\alpha$-multistable marginal Lévy measures determine the local jump behavior of the individual degradation components, whereas the Clayton Lévy copula determines the dependence among component-wise jumps. Although the Lévy-copula specification is time-invariant, the resulting joint local Lévy measure is time-dependent because its marginal Lévy measures vary over time $s$.

The cumulative joint Lévy measure over the observation interval [0, $t$] is therefore obtained by integrating the joint local Lévy measure over time:

$$v_t^C(dy) = \int_0^t v_s^C(dy)ds. \tag{25}$$

This cumulative measure determines the finite-time joint distribution of the degradation process through the Lévy–Khintchine representation developed in Section 3.2.

### 4.4. Marginal and System Reliability under Clayton Distributional Dependence

Section 3.3 establishes the general first-passage reliability formulation for multidimensional nonhomogeneous subordinators. For the $k$th $\alpha$-multistable marginal degradation process considered here, the marginal reliability is given by Eq. (16):$R_k(t) = F_k(x_k, t)$. The marginal cumulative distribution function appearing in Eq. (16) is defined in Eq. (14). It is derived from the finite-time marginal characteristic function in Eq. (11) using the Fourier inversion formula in Eq. (13).

To obtain an explicit system reliability expression, we adopt a Clayton distributional copula for the finite-time degradation variables. Let $C_F$ denote the distributional copula, which is distinct from the Lévy copula $C_L$ introduced in Section 4.2. For $v_k \in [0,1]$, the $m$-dimensional Clayton distributional copula is

$$C_F(v_1, \dots, v_m\,;\, \theta_F) = \left(\sum_{k=1}^{m} v_k^{-\theta_F} - m + 1\right)^{\frac{-1}{\theta_F}}, \quad \theta_F > 0. \tag{26}$$

Substituting $v_k = R_k(t)$ into the general reliability representation in Section 3.3 gives:

$$R(t) = C_F(R_1(t), \dots, R_m(t); \theta_F) = \left(\sum_{k=1}^{m} R_k(t)^{-\theta_F} - m + 1\right)^{\frac{-1}{\theta_F}}. \tag{27}$$

Equation (27) gives the system reliability of the $\alpha$-multistable degradation model under Clayton distributional dependence. The marginal reliability functions $R_k(t)$ describe the time-dependent degradation behavior of the individual components, whereas $\theta_F$ governs dependence among the finite-time degradation variables. The formulation therefore distinguishes marginal nonhomogeneous degradation from system-level dependence.

Before proceeding to the numerical analysis, we verify two basic consistency properties of the proposed reliability formulation. The first concerns the initial condition of the degradation process, while the second confirms that the homogeneous model is recovered as a special case when the stability parameters are constant over time.

**Initial-condition consistency.** Since every component degradation process starts from zero, no degradation has accumulated at time zero. Consequently, all marginal reliability functions equal one at the initial time. Substituting these marginal reliabilities into Eq. (27) gives a system reliability of one, which is consistent with the initial condition of the degradation process.

**Homogeneous limiting case.** If the time-dependent stability parameter of each component is constant over time, the marginal Lévy measures become time-invariant, and the proposed nonhomogeneous formulation reduces to the corresponding homogeneous positive α-stable degradation model. Therefore, the homogeneous model is a special case of the framework developed in this section.

## 5. Simulation and Numerical Examples

This section applies the theoretical framework developed in Sections 3 and 4 to a two-dimensional $\alpha$-multistable degradation system. The numerical study examines three features of the proposed framework: (*i*) the effect of the Lévy-copula dependence parameter on the joint evolution of the degradation trajectories, (*ii*) the time-varying jump behavior generated by the nonhomogeneous marginal Lévy measures, and (*iii*) the effects of the marginal degradation distributions and finite-time dependence on system reliability. All model definitions and analytical derivations are presented in Sections 3 and 4. This section focuses on the numerical implementation, simulation results, and interpretation of the resulting degradation and reliability behavior.

### 5.1. Simulation Setup

The numerical study considers a two-dimensional nonhomogeneous degradation process consisting of two α-multistable marginal subordinators. The marginal local Lévy measures are specified according to Eq. (21), and the joint local Lévy measure is evaluated by applying the two-dimensional version of Eq. (24) with the Clayton Lévy copula defined in Eq. (22).

To illustrate the effect of jump dependence, the Clayton Lévy-copula parameter $\theta^L$ is assigned three values: $\theta = 15$, $\theta = 12$ and $\theta = 5$. The marginal model parameters are held fixed across the three settings. Figures 1–3 present representative simulated realizations under the three dependence settings. Because each figure displays a single realization, these trajectories are intended to provide a qualitative illustration of dependence behavior rather than a formal quantitative comparison.

### 5.2. Effect of Lévy-Copula Dependence on Degradation Paths

Figures 1-3 present representative two-dimensional degradation trajectories generated under $\theta_L = 15, 12$, and $5$, respectively. In all three cases, the degradation trajectories are non-decreasing and exhibit nonuniform jump behavior over time, consistent with the $\alpha$-multistable marginal specification.

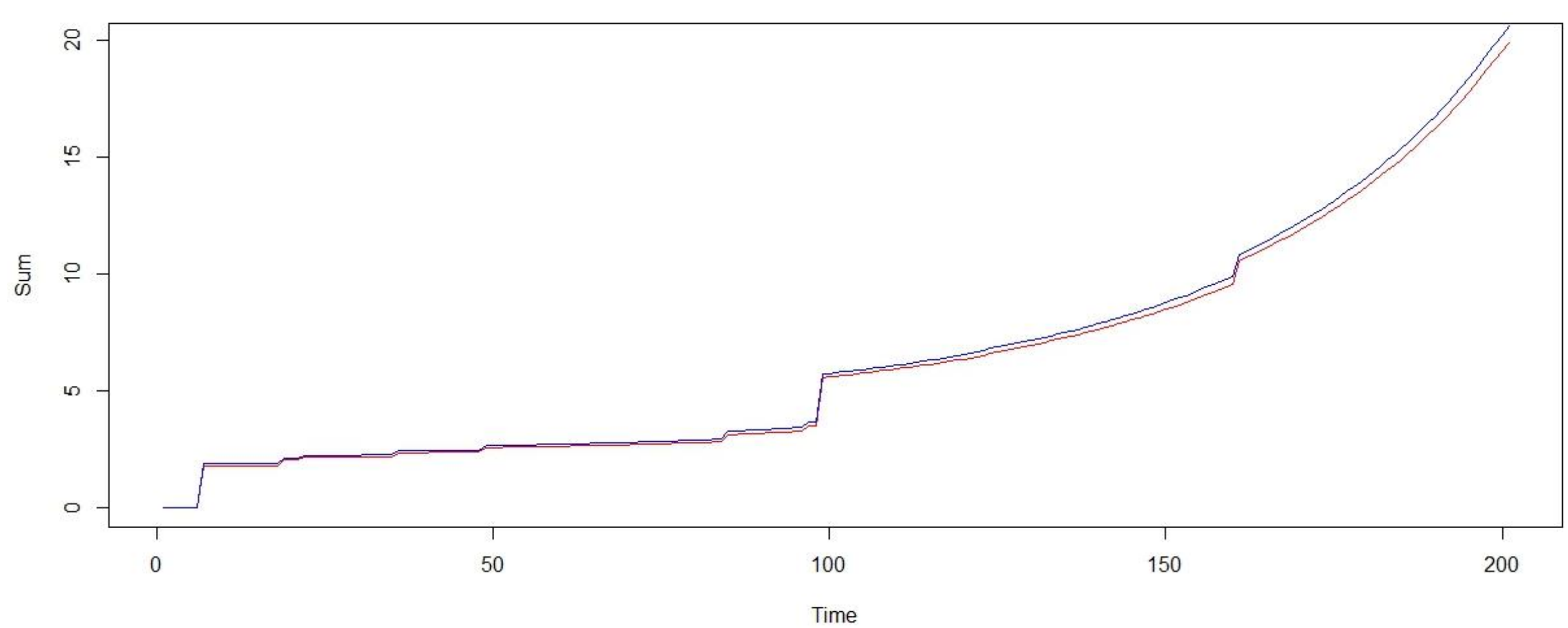

**Figure 1**: Two-dimensional nonhomogeneous subordinator ($\theta_L = 15$)

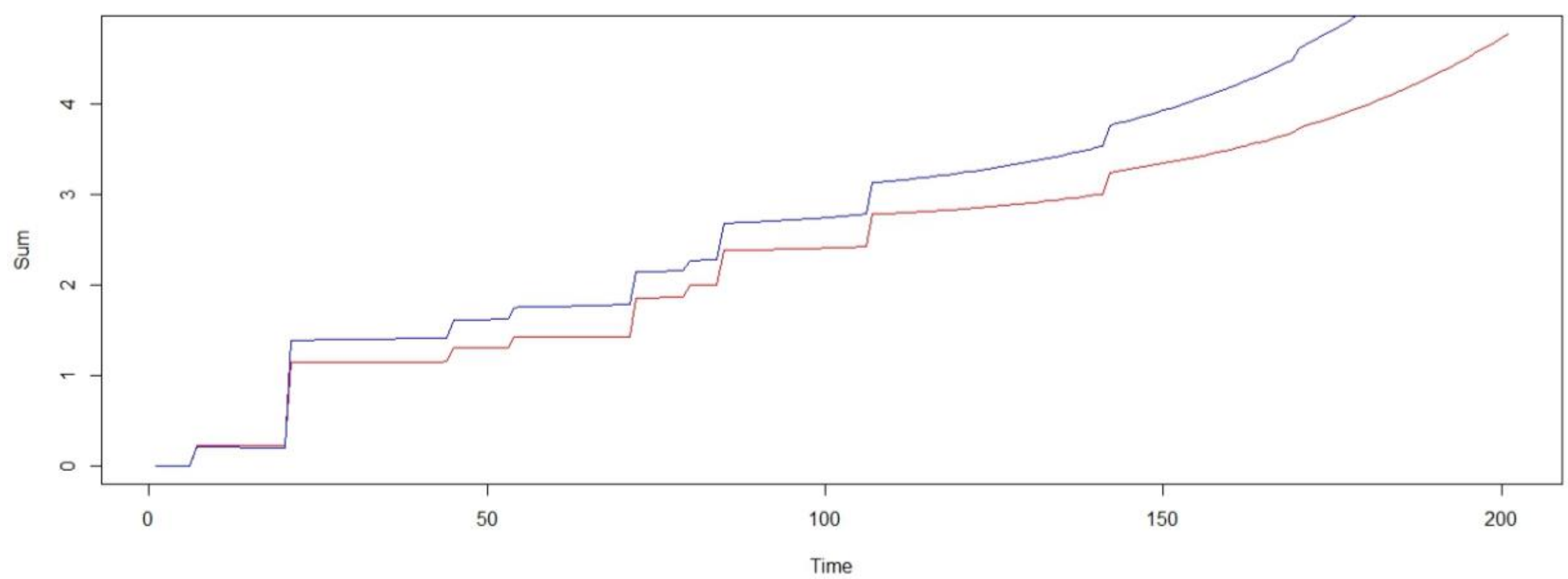

**Figure 2**: Two-dimensional nonhomogeneous subordinator ($\theta_L = 12$)

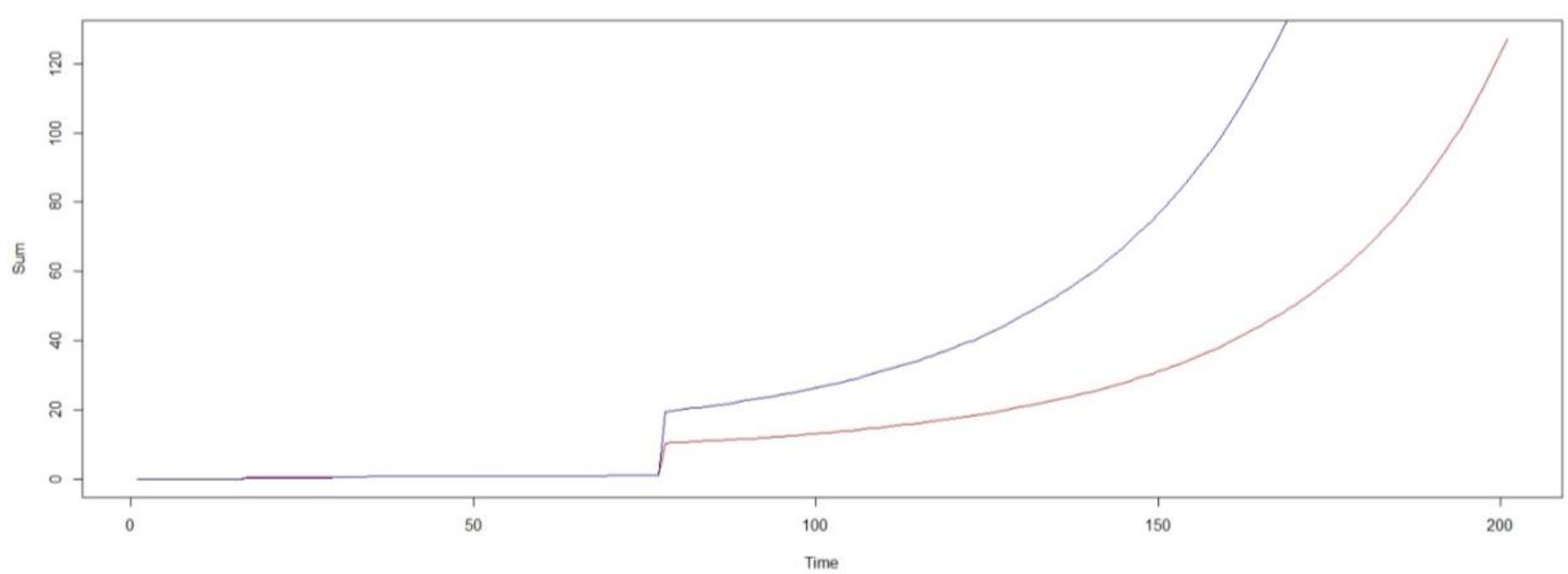

**Figure 3**: Two-dimensional nonhomogeneous subordinator ($\theta_L = 5$)

Figures 1–3 provide qualitative illustrations of degradation trajectories under different values of $\theta^L$. In the representative realizations shown here, the trajectories exhibit different degrees of synchronization and separation between component-wise jumps. The $\alpha$-multistable marginal specification governs the time-varying jump behavior of each component, whereas the Clayton Lévy copula governs dependence between their jumps. Because each panel presents a single

simulated realization, differences in the cumulative degradation magnitudes should not be attributed solely to the dependence parameter.

### 5.3. Time-Varying Marginal Lévy Measure

In addition to cross-component dependence, the nonhomogeneous behavior of the degradation process is governed by the time-varying marginal Lévy measures. To isolate this feature, we examine the local marginal Lévy measure at a fixed degradation jump magnitude while varying the time $s$.

Figure 4 shows that the marginal local Lévy measure changes substantially over time. Under the parameter specification used in this numerical example, the local jump intensity initially increases, reaches a maximum near $s = 1.2$, and subsequently decreases. Thus, the jump characteristics of the degradation process are not constant across the observation period.

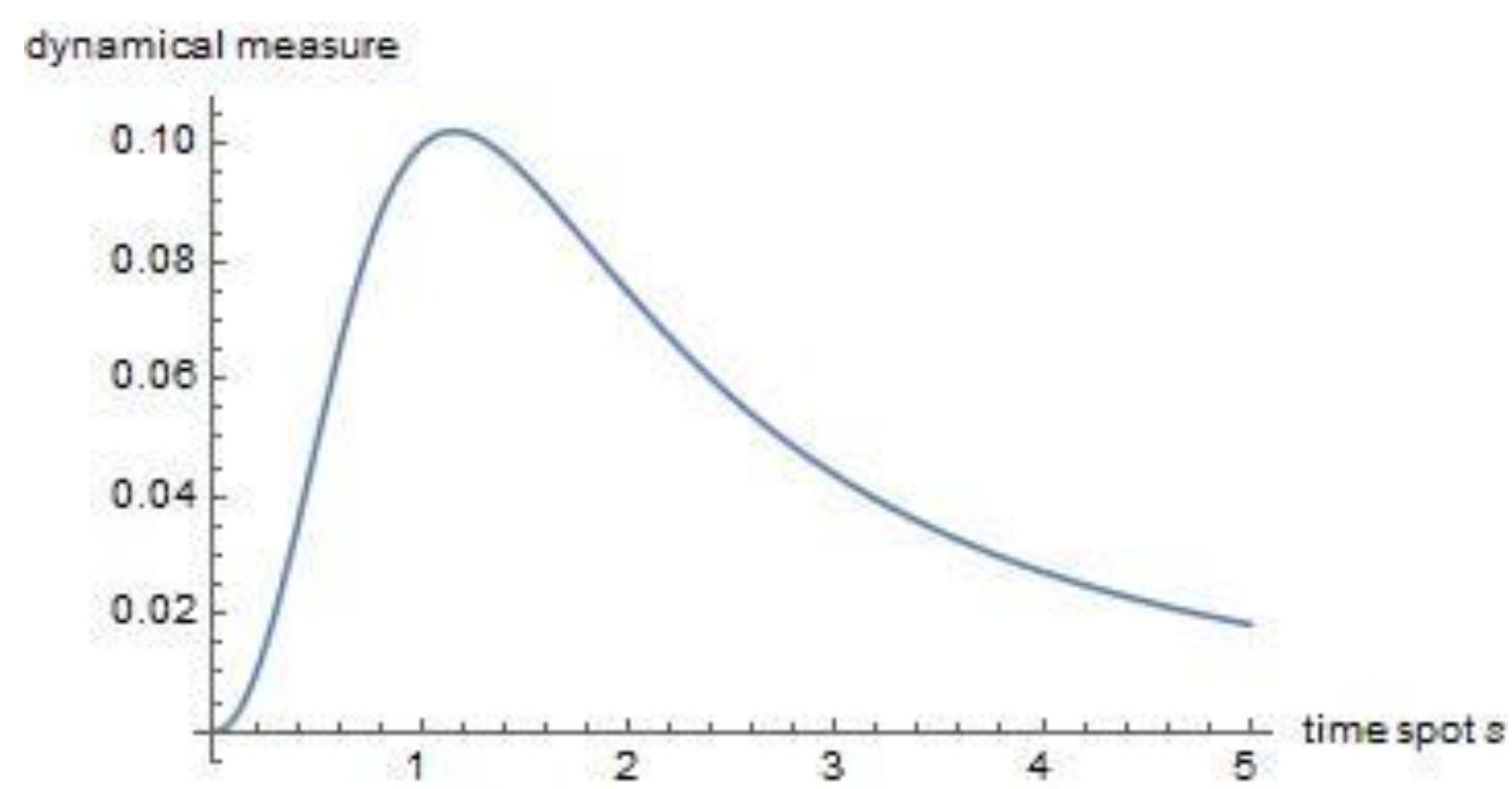

**Figure 4**: Time-varying marginal dynamic Lévy measure

This behavior provides a direct numerical illustration of the nonhomogeneous structure introduced in Section 4.1. In contrast, for a homogeneous Lévy subordinator, the local Lévy measure would not vary with time. Therefore, Figure 4 demonstrates how the $\alpha$-multistable specification allows the local degradation dynamics to evolve over time.

### 5.4. Marginal and System Reliability

To illustrate how the time-dependent degradation dynamics are translated into reliability measures, we evaluate the marginal reliability of the two components and the corresponding system reliability. The marginal reliability functions $R_1(t)$ and $R_2(t)$ are calculated from the finite-time marginal degradation distributions obtained using the formulations developed in Sections 3.2 and 4.1.

For the system-level reliability calculation, we use the Clayton distributional copula introduced in Section 4.4 with the finite-time dependence parameter set to $\theta_F = 5$ . The corresponding system reliability is then calculated using Eq. (27).

Figure 5 presents the resulting system reliability function, while Figure 6 shows the two marginal reliability functions. Both the marginal and system reliability functions decrease monotonically with observation time as cumulative degradation progresses. This monotonic decline is consistent with the non-decreasing sample paths of the underlying α-multistable subordinators.

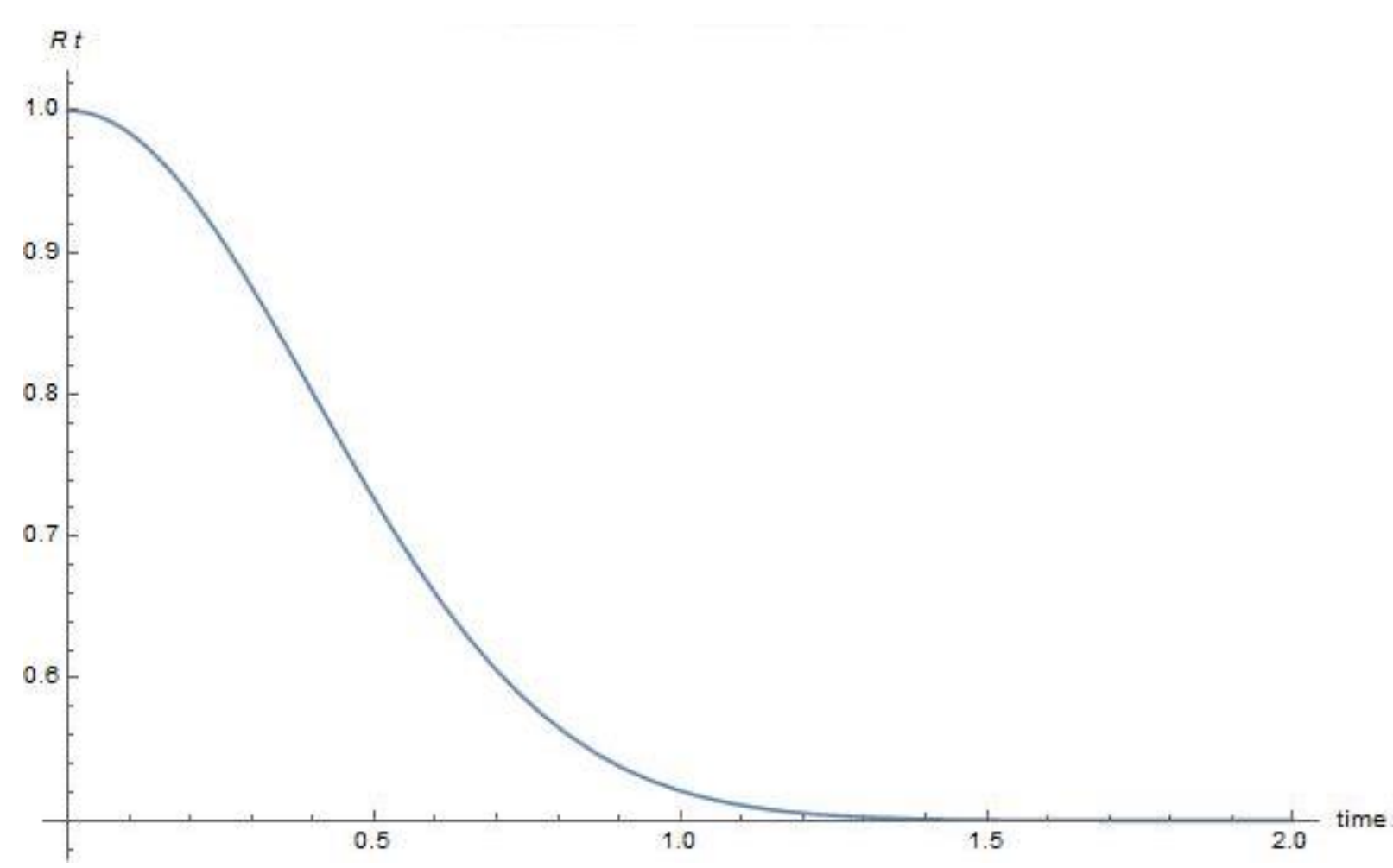

**Figure 5:** System reliability under Clayton distributional dependence with $\theta^F = 5$.

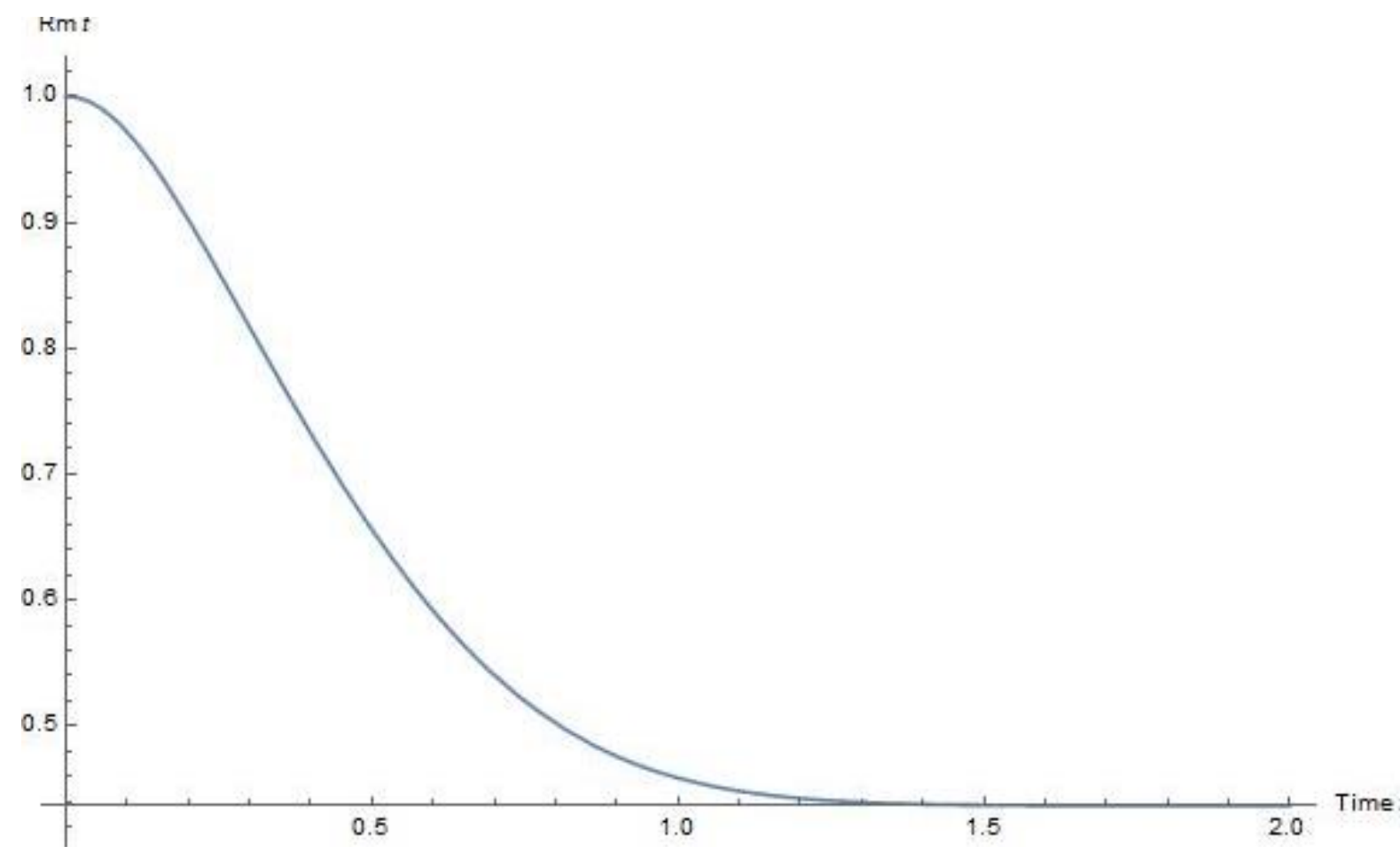

**Figure 6:** Marginal reliability functions of the two degradation components.

The system reliability incorporates two distinct sources of information. The marginal reliability functions reflect the time-dependent degradation dynamics of the individual components, whereas the Clayton distributional copula accounts for dependence among the finite-time degradation variables. Because the degradation system is assumed to operate as a series system, system survival requires all degradation components to remain below their prescribed thresholds. Therefore, the system reliability provides a joint measure that combines marginal degradation behavior and finite-time dependence.

Together with the time-varying Lévy measure illustrated in Figure 4, these results demonstrate how local nonhomogeneous jump characteristics propagate into finite-time degradation distributions and ultimately affect system-level reliability.

## 6. Conclusion and Discussion

This study develops a multi-dimensional reliability framework for degradation processes with time-dependent random jumps by extending homogeneous Lévy-subordinator models to the nonhomogeneous setting. Each marginal degradation process is modeled by a nonhomogeneous subordinator characterized by a time-dependent Lévy measure, while dependence among component-wise jumps is represented through a Lévy-copula structure. By allowing the local Lévy characteristics to evolve over time, the proposed formulation relaxes the stationary-increment assumption commonly adopted in conventional Lévy degradation models.

A general distributional framework is first established for multi-dimensional nonhomogeneous subordinators. The marginal and joint local Lévy characteristics are connected through a time-dependent Lévy-copula representation, and their accumulation over the observation interval determines the finite-time degradation distributions through the Lévy–Khintchine framework. The corresponding marginal characteristic functions, probability density functions, and cumulative distribution functions are then used to formulate the reliability for a multi-dimensional series degradation system based on first-passage times. At the finite-time distributional level, system reliability is represented through a distributional copula, which is distinct from the Lévy copula governing jump-level dependence.

The general framework is further specialized to $\alpha$-multistable marginal subordinators with Clayton dependence. In this specialization, the time-varying $\alpha$-multistable Lévy measures characterize the non-stationary marginal degradation behavior, while the Clayton Lévy copula describes dependence among component-wise jumps. A Clayton distributional copula is then used to represent dependence among the finite-time degradation variables in the system reliability formulation. The homogeneous positive $\alpha$-stable degradation model is recovered as a special case when the time-varying stability parameters are constant.

The numerical examples illustrate the main features of the proposed framework. The simulated degradation paths are non-decreasing and exhibit time-varying jump behavior, consistent with the $\alpha$-multistable specification. The representative simulations also illustrate how changes in Lévy-copula dependence strength affect the degree of co-movement between degradation components. The marginal local Lévy measure varies over time, directly demonstrating the nonhomogeneous jump mechanism. Finally, both the marginal and system reliability functions decrease as cumulative degradation progresses, illustrating how time-varying degradation dynamics propagate into finite-time reliability measures.

The proposed framework provides a flexible basis for reliability and lifetime analysis of systems involving multiple dependent degradation processes. Several extensions remain for future research. For high-dimensional systems, dimension-reduction or feature-selection methods may be incorporated to reduce model complexity before constructing high-dimensional dependence models. Future work may also consider time-varying Lévy-copula parameters, alternative families of Lévy and distributional copulas, and data-driven estimation of time-dependent marginal Lévy measures from observed degradation data.

## References

[1] Y. Shu, Q. Feng, and D. W. Coit (2015) Life Distribution Analysis Based on Lévy Subordinators for Degradation with Random Jumps, *Naval Research Logistics*, **(62)** 483-492.

[2] Y. Shu, Q. Feng, E. Kao and H. Liu (2016) Lévy-driven Non-Gaussian Ornstein-Uhlenbeck Processes for Degradation-based Reliability Analysis, *IIE Transactions* **48(11)** 993-1003.

[3] Y. Shu, Q. Feng, E. Kao, D. W. Coit and H. Liu (2021) Markov additive processes for degradation with jumps under dynamic environments, *Naval Research Logistics* **(68)** 908-919.

[4] Y. Shi, Q. Feng, Y. Shu and Y. Xiang (2020) Multi-dimensional Lévy processes with Lévy copulas for multiple dependent degradation processes in lifetime analysis, *Quality Engineering*, **32(3)** 434-448.

[5] Y. Shi and Q. Feng (2024) Reliability analysis for binary Lévy degradation processes using characteristic functions, *Quali Reliab Engng Int*, **2024** 1-18.

[6] L. MacLean, A. Richman, and M. Hudak (2018) Failure Rates for Aging Aircraft, *Safety* 4(1):7. https://doi.org/10.3390/safety4010007.

[7] D. M. Boyajian (2002) Mode I fracture and durability of the CFRP -concrete interface bond, *Graduate Theses, Dissertations, and Problem Reports*, **2459** https://researchrepository.wvu.edu/etd/2459

[8] D. Liu, S. Wang and M. Tomovic (2020) Degradation modeling method for rotary lip seal based on failure mechanism analysis and stochastic process, *Eksploatacja i Niezawodność – Maintenance and Reliability*, 22(3), pp.381-390. https://doi.org/10.17531/ein.2020.3.1

[9] R. Cont and P. Tankov (2004) *Financial Modeling with Jump Processes*, Chapman & Hall/CRC.

[10] F. Black and M. Scholes (1973) The Pricing of Options and Corporate Liabilities, *the Journal of Political Economy*, **81(3)** 637-654.

[11] M. Abdel-Hameed (2014) *Lévy processes and their applications in reliability and storage*, Springer Briefs in Statistics.

[12] E. Orsingher, C. Ricciuti and B. Toaldo (2016) Time-inhomogeneous jump processes and variable order operators, *Potential Analysis*, published online.

[13] L. Beghin and C. Ricciuti (2018) Time-inhomogeneous fractional Poisson processes defined by the multistable subordinator, *Mathematics*, DOI:10.1080/07362994.2018.1548970.

[14] K. Sato (1999) *Lévy processes and Infinitely Divisible Distributions*, Cambridge University Press.

[15] K. Aas, C. Czado, A. Frigessi and H. Bakken (2009) Pair-copula constructions of multiple dependence, *Insurance: Mathematics and Economics* **44** 182–198.

[16] N.G. Shephard (1991) From characteristic function to distribution function: A simple framework for the theory, *Econ Theory* **7** 519-529.

[17] W. Kluge (2005) Time-inhomogeneous Levy processes in interest rate and credit risk models, *Dissertation zur Erlangung des Doktorgrades der Fakult¨at f¨ur Mathematik und Physik der Albert-Ludwigs-Universit¨at Freiburg im Breisgau*.

[18] E. Cuvelier, M. Noirhomme-Fraiture (2005) Clayton copula and mixture decomposition, In J. Janssen, & P. Lenca (Eds.), *Applied Stochastic Models and Data Analysis (ASMDA 2005)*, Brest 17-20.

[19] T. Castro and L. Landesa (2024) A dependent complex degrading system with non-periodic inspection times, https://arxiv.org/abs/2401.09141

[20] X. Sun, W. Cai and M. Li (2021) A hierarchical modeling approach for degradation data with mixed-type covariates and latent heterogeneity, *Reliability Engineering & Systems Safety*, **216** 0951-8320.

[21] N. D. Singpurwalla (1995) Survival in Dynamic Environments, *Statistical Science*, **10(1)** 86-103.

[22] Mozumder, M. S. U. (2011) Some Derivations on Levy and Jump Processes. GANIT: Journal of Bangladesh Mathematical Society, 29, 11–22. https://doi.org/10.3329/ganit.v29i0.851.

[23] William E. Pruitt (1981) The Growth of Random Walks and Levy Processes, *The Annals of Probability*, Ann. Probab. **9(6)** 948-956.